\documentclass[11pt]{amsart}
\usepackage{amsmath}
\usepackage{amssymb}
\usepackage{amsfonts}
\usepackage{graphicx}

\DeclareSymbolFont{AMSb}{U}{msb}{m}{n}
\DeclareMathSymbol{\N}{\mathbin}{AMSb}{"4E}
\DeclareMathSymbol{\Z}{\mathbin}{AMSb}{"5A}
\DeclareMathSymbol{\R}{\mathbin}{AMSb}{"52}
\DeclareMathSymbol{\Q}{\mathbin}{AMSb}{"51}
\DeclareMathSymbol{\I}{\mathbin}{AMSb}{"49}
\DeclareMathSymbol{\C}{\mathbin}{AMSb}{"43}

\theoremstyle{definition}

\theoremstyle{corollary}

\theoremstyle{example}

\theoremstyle{note}

\theoremstyle{notation}

\numberwithin{equation}{section}
\begin{document}
\title{DYNAMICS NEAR AN IDEMPOTENT}

\author{Md Moid Shaikh and Sourav Kanti Patra}
\address{Md Moid Shaikh, Department of Mathematics, Maharaja Manindra Chandra College,
20, Ramkanto Bose Street, Kolkata-700 003, West Bengal, India}
\email{mdmoidshaikh@gmail.com}

\address{Sourav Kanti Patra, Department of Mathematics, Ramakrishna Mission Vidyamandira,
Belur Math, Howrah-711202, West Bengal, India}
\email{souravkantipatra@gmail.com}


\begin {abstract}
Hindman and Leader first introduced the notion of semigroup of ultrafilters converging to zero for 
a dense subsemigroups of $((0,\infty),+)$. Using the algebraic structure of the Stone-$\breve{C}$ech 
compactification, Tootkabani and Vahed generalized and extended this notion to an idempotent instead 
of zero, that is a semigroup of ultrafilters converging to an idempotent $e$ for a dense subsemigroups 
of a semitopological semigroup $(T, +)$ and they gave the combinatorial proof of central set theorem near 
$e$. Algebraically one can also define quasi-central sets near $e$ for dense subsemigroups of $(T, +)$. 
In a dense subsemigroup of $(T,+)$, C-sets near $e$ are the sets, which satisfy the conclusions of the 
central sets theorem near $e$. S. K. Patra gave dynamical characterizations of these combinatorially rich 
sets near zero. In this paper we shall prove  these dynamical characterizations for these combinatorially 
rich sets near $e$.

AMS subjclass [2010] : 37B20; 37B05; 05B10.
\end{abstract}

\maketitle

\section {introduction}
In topological dynamics there are many notions of large subsets of a semigroup are defined. Many of them are 
combinatorially rich and playing a significant role in Ramsey Theory. One of them is central set which was 
first introduced by H. Furstenberg using the notions of proximal and uniformly recurrent points from topological 
dynamics in \cite[Definition 8.3]{frus}. Bergelson and Hindman later provided an algebraic characterization of 
central sets of $\mathbb{N}$ in \cite[Section 6]{berg}.
This algebraic characterization of central sets is nicely extended to any semigroup in terms of the algebraic
structure of $\beta S$, where $\beta S$ is the Stone-$\breve{C}$ech compactification of the discrete semigroup
$(S,\cdot)$. Before giving this characterization, we present a brief description of the algebraic structure of 
$\beta S$ for a discrete semigroup $(S,\cdot)$. For details one can see \cite{hindalg}.\\

Let $(S,\cdot)$ be any discrete semigroup then the Stone-$\breve{C}$ech compactification $\beta S$ of the discrete 
semigroup $(S,\cdot)$ is the set of all ultrafilters on $S$, the principal ultrafilters being identified with the 
points of $S$. Given $A\subseteq S$ let us set, $\bar{A}=\{ p\in \beta S : A\in p\}$. Then the set $\{ \bar{A} : 
A\subseteq S\}$ is a  clopen basis for a topology on $\beta S$. The operation $\cdot$ on $S$ can be extended to the 
Stone-$\breve{C}$ech compactification $\beta S$ of $S$ so that $(\beta S,\cdot)$ is a compact right topological 
semigroup (meaning that for any $p\in \beta S$, the function $\rho_p:\beta S \rightarrow \beta S$ defined by 
$\rho_p(q)=q\cdot p$ is continuous) with $S$ contained in its topological center (meaning that for any $x\in S$ 
the function $\lambda_x : \beta S \rightarrow \beta S$ 
defined by $\lambda_x (q)=x\cdot q$ is continuous). Given $p,q \in \beta S $ and 
$A \in S$, $A \in p\cdot q$ if and only if $\{ x \in S :x^{-1}A \in q\} \in p$, where 
$x^{-1}A=\{ y \in S: x\cdot y \in A \}$. A non-empty subset $I$ of a semigroup $(T,\cdot)$ 
is called a left ideal of $S$ if $T\cdot I \subseteq I$, a right ideal of $S$ if 
$I\cdot T \subseteq I$ and a two sided ideal (or simply an ideal) if it is both a 
left and a right ideal. A minimal left ideal is a left ideal that does not 
contain any proper left ideal. Similarly we can define minimal right ideal 
and the smallest ideal. Any compact Hausdorff right topological semigroup $(T,\cdot)$ 
has the smallest two sided ideal\\
$$\begin{array}{ccc}
K(T) & = & \bigcup\{L:L \text{ is a minimal left ideal of } T\} \\
& = & \,\,\,\,\,\bigcup\{R:R \text{ is a minimal right ideal of } T\}\\
\end{array}$$.\\
Given a minimal left ideal $L$ and a minimal right ideal $R$, $L\cap R$ is a group, and in particular contains an idempotent.
An idempotent is minimal if and only if it is a member of the smallest ideal.\\ 
Now we present Bergelson's characterizations of Central sets.\\

\textbf{Definition 1.1.} Let $S$ be a discrete semigroup and let $C$ be a subset of $S$. 
Then $C$ is central if and only if there is an idempotent $p$ in $K(\beta S)$ such 
that $C \in p$.\\
We now state the well-known definition of dynamical system.\\

\textbf{Definition 1.2.} A dynamical system is a pair $(X,\langle T_s\rangle _{s \in S })$ 
such that

(i) $X$ is compact Hausdorff space,

(ii) $S$ is a semigroup,

(iii) for each $s\in S$, $T_s:X\rightarrow X$ and $T_s$ is continuous, and

(iv) for all $s,t$ , $T_s\circ T_t=T_{st}$.\\

 We write $P_f(X)$ is the set of finite non-empty subsets of $X$, for any set $X$.\\

The following definitions are very essential to give a dynamical characterization of central sets in an arbitrary semigroup $S$ .\\

\textbf{Definition 1.3.} (\cite[Definition 3.1]{hindcen}) Let $S$ be a semigroup and 
let $A\subseteq S$.

(a) The set $A$ is syndetic if and only if there is some $G \in P_f(S)$ such that 
    $S=\bigcup_{t \in G}t^{-1}A $.

(b) The set $A$ is piece-wise syndetic if and only if there is some  $G \in P_f(S)$ 
such that for any $F \in P_f(S)$ there is some $x \in S$ with 
$Fx\subseteq \bigcup_{t \in G}t^{-1}A $.\\

Recall the definitions of proximality and uniform recurrence in a dynamical system 
from \cite[Definition 1.2(b)]{burns} and \cite[Definition 1.2(c)]{burns}.\\

\textbf{Definition 1.4.} Let $(X,\langle T_s\rangle _{s \in S})$ is a dynamical system.

(a) A point $y \in S$ is uniformly recurrent if and only if for every neighbourhood 
    $U$ of $y$, $\{ s\in S : T_s(y) \in U  \}$ is syndetic.

(b) The points $x$ and $y$ of $X$ are proximal if and only if for every neighbourhood 
    $U$ of the diagonal in $X\times X$, there is some $s \in S$ such that 
    $(T_s(x), T_s(y)) \in U$.\\
    
By \cite[Theorem 2.4]{shi}, a subset $C$ of a semigroup $S$ is central if and only if 
there exist a dynamical system $(X,\langle T_s\rangle _{s \in S})$, points $x$ and $y$ of $X$ 
and a neighbourhood $U$ of $y$ such that $y$ is uniformly recurrent, $x$ and $y$ 
are proximal and $C=\{s  \in S: T_s(x) \in U  \}$.\\

Now we discuss some basic definitions, conventions and results for dynamical 
characterization of members of certain idempotent ultrafilters.\\

We state \cite[Definition 2.1]{john}.\\

\textbf{Definition 1.5.}  Let $S$ be a nonempty discrete space and let
$\mathcal{K}$ be a filter on $S$.

(a) $\bar{\mathcal{K}}= \{p \in \beta S : \mathcal{K} \subseteq p \}$.

(b) $\mathcal{L}(\mathcal{K})=\{A\subseteq S: S\setminus A \not \in \mathcal{K}\}$.\\

In \cite[Theorem 3.20]{hindcen} it is proved that the function $\mathcal{K} \rightarrow \bar{\mathcal{K}}$ 
is a bijection from the collection of all filters on $S$ onto the collection of all compact subspaces of $\beta S$.

 Relating the above two concepts we have the following important theorem.\\
 
\textbf{Theorem 1.6.}(\cite[Theorem 1.6]{sou}) Let $S$ be a nonempty discrete space and $\mathcal{K}$ a filter on $S$.
(a) $\bar{\mathcal{K}}=\{ p \in \beta S: A \in \mathcal{L}(\mathcal{K})$ for all $A \in p$\}. 

(b) Let $\mathcal{B} \subseteq \mathcal{L}(\mathcal{K})$ be closed under finite intersections. Then 
there exists a $p \in \beta S$ with $\mathcal{B} \subseteq p \subseteq \mathcal{L}(\mathcal{K})$.\\

\textbf{Proof.} Both of these assertions follow from \cite[Theorem 3.11]{hindalg}.\\
We recall (\cite[Definition 3.1]{john}).\\

\textbf{Definition 1.7.} Let $(X,\langle T_s\rangle _{s \in S})$ be a 
dynamical system, $x$ and $y$ points in $X$, and $\mathcal{K}$ a filter on $S$. 
The pair $(x,y)$ is called jointly $\mathcal{K}$-recurrent if and only if for 
every neighbourhood $U$ of $y$ we have $\{ s \in S: T_s(x) \in U$ and 
$T_s(y) \in U \} \in \mathcal{L}(\mathcal{K})$.\\

Relating to above definition the following theorem  was obtained independently (for filters on $\mathbb{N}$) in 
\cite[Theorem 4.11]{jli} and (again independently) for filters on an arbitrary semigroup
in \cite[Theorem 5.2.3]{chris}.\\

\textbf{Theorem 1.8.} Let $S$ be a semigroup, $\mathcal{K}$ be a filter on 
$S$ such that $\bar{\mathcal{K}}$ is a compact subsemigroup of $\beta S$, and let 
$A \subseteq S$. Then $A$ is a member of an idempotent in $\bar{\mathcal{K}}$ if 
and only if there exists a dynamical system $(X,\langle T_s\rangle _{s \in S})$ with points 
$x$ and $y$ in $X$ and there exists a neighbourhood $U$ of $y$ such that the 
pair $(x,y)$ is jointly $\mathcal{K}$-recurrent and $A=\{ s \in S: T_s(x) \in U \}$.\\

\textbf{Proof.} See the proof of \cite[Theorem 3.3]{john}.\\ 

Now we define another combinatorially rich sets, quasi-central sets.\\

\textbf{Definition 1.9.} (\cite[Definition 1.2]{hindcen}) Let S be a discrete semigroup and let 
$C$ be a subset of $S$. Then $C$ is quasi-central if and only if there is an idempotent 
$p$ in $ClK(\beta S)$ such that $C \in p$.\\

Now recall \cite[Definition 4.1 and Definition 4.4]{john}.\\

\textbf{Definition 1.10.} Let $S$ be a semigroup.

(a) For each positive integer $m$ put $J_m=\{ (t_1,t_2,...,t_m)\in 
\mathbb{N}^m:t_1<t_2<...<t_m \}$.

(b) Given $m \in \mathbb{N}, a \in S^{m+1}, t \in J_m$, and $t \in \tau $, put 
$x(m,a,t,f)=\prod^m_{i=1}(a(i)f(t_i))a(m+1) $ where $\tau=S^{\mathbb{N}}$.

(c) We call a subset $A\subseteq S$, a $C$-set if and only if there exist 
functions $m:P_f(\tau)\rightarrow \mathbb{N}$, $\alpha \in \prod_{{F \in P_f(\tau)}}{S^{m(F)+1}}$, 
and $\tau \in \prod_{{F \in P_f(\tau)}}{J_{m(F)}}$ such that the following two 
statements are satisfied:

(1) If $F,G \in P_f(\tau)$ and $F\subsetneq G$ then $\tau(F)(m(F))<\tau(G)(1)$.

(2) Whenever $m \in \mathbb{N}, G_1, G_2,...,G_m$ is a finite sequence in 
$P_f(\tau)$ with $G_1 \subsetneq G_2 \subsetneq  ...\subsetneq G_m$ and for 
each $i \in \{ 1,2,...,m \}, f_i \in G_i$ then we have\\
$\prod^m_{i=1}x(m(G_i), \alpha (G_i), \tau (G_i),f_i)\in A$.

(d) We call a subset $A\subseteq S$, a $J$-set if and only if for every 
$F \in P_f(\tau)$, there exist $m \in \mathbb{N},a \in S^{m+1}$ and 
$t \in J_m$ such that for all $f \in F$, $x(m,a,t,f)\in A$.\\

The following definitions are very useful to give dynamical characterizations of quasi-central sets 
and $C$-sets.\\

\textbf{Definition 1.11.}(\cite[Definition 1.11]{sou}) Let $(X,\langle T_s\rangle _{s \in S})$ be a dynamical system and 
let $x,y \in X$.

(a) The pair $(x,y)$ is jointly intermittently uniformly recurrent 
(abbreviated as $JIUR$) if and only if for every neighbourhood $U$ of $y$, 
$\{ s \in S: T_s(x)\in U$ and $T_s(y)\in U \}$ is piecewise syndetic.

(b) The pair $(x,y)$ is jointly intermittently almost uniform recurrent 
(abbreviated as $JIAUR$) if and only if for every neighbourhood $U$ of $y$, 
$\{ s \in S: T_s(x)\in U$ and $T_s(y)\in U \}$ is a $J$-set.\\

Using Theorem 1.8, one can get dynamical characterizations of quasi-central sets 
and $C$-sets in terms of $JIUR$ and $JIAUR$ respectively.\\

\textbf{Theorem 1.12.} Let S be a semigroup and let $C\subseteq S$. The set $C$ 
is quasi-central if and only if there exist a dynamical system 
$(X,\langle T_s\rangle _{s \in S})$, points $x$ and $y$ in $X$, and a neighbourhood 
$U$ of $y$ such that the pair $(x,y)$ is $JIUR$ and $C=\{ s \in S: T_s(x)\in U\}$.\\

\textbf{Proof.} See the proof of \cite[Theorem 3.4]{burns}.\\

\textbf{Theorem 1.13.} Let S be a semigroup and let $C\subseteq S$. The set $C$ 
is a $C$-set if and only if there exist a dynamical system $(X,\langle T_s\rangle _{s \in S})$, 
points $x$ and $y$ in $X$, and a neighbourhood $U$ of $y$ such that the 
pair $(x,y)$ is  $JIAUR$ and $C=\{ s \in S: T_s(x)\in U  \}$.\\

\textbf{Proof.} This is established by the proof of \cite[Theorem 4.8]{john}.\\

 In \cite{sou} Patra studied the above results in near zero and in this paper we shall extend these cases near an idempotent.\\ 
In \cite{hindult} Hindman and Leader discuss the semigroup of ultrafilters converging to zero
for a dense subsemigroups of $((0,\infty),+)$ in details.\\
Now we present briefly the semigroup of ultrafilters converging to an idempotent $e$ for a dense subsemigroups of a semitopological semigroup
$(T, +)$. See \cite{mat} for details.\\

Let $(T, +)$ be a Hausdorff semitopological semigroup. Let $E(T)$ denote the set of all idempotents in $T$.
For every $x \in T$, $\tau_x$ is the collection of all neighbourhoods of $x$.

Now we shall consider semigroups which are dense in $(T,+)$ with usual topology.\\

\textbf{Definition: 1.14.}(\cite[Definition 2.1]{mat})Let $(T, +)$ be a semitopological semigroup and $S$ is a dense subsemigroup of
$T$.\\
(a) Given $x\in T$, $x^{*}_S=\{ p \in \beta S_d : x\in \bigcap_{A\in p} cl_T(A)\}$.\\
(b) $B(S)=\bigcup _{x\in T}x^{*}$.\\
(c) ${\infty}^{*}= \beta S_d\backslash B(S)$.\\

\textbf{Definition: 1.15.} If $S$ is a dense subsemigroup of $(T,+)$\\ 
then $e^{*}_S=\{ p \in \beta T_d : e\in \bigcap_{A\in p} cl_T(A)\}$. ($T_d$ represents the set $T$ with discrete topology.)\\

It was proved in \cite[Lemma 2.3]{mat} that $e^{*}_S$ is a compact right topological
subsemigroup of $(\beta S_d,+)$. Using \cite[Lemma 2.2]{mat} we can say that $e^{*}_S$ is disjoint 
from $K(\beta S_d)$ and hence gives some new information which is not available 
from $K(\beta S_d)$.\\

As $e^{*}_S$ is a compact right topological semigroup, $e^{*}_S$ has the smallest ideal 
$K=K(e^{*}_S)$ and hence $e^{*}_S$ contains an idempotent \cite[Corollary 2.10]{elli}. One can define minimal idempotent also in $e^{*}_S$.\\

Therefore one can have combinatorially rich sets near an idempotent such as central sets, quasi-central sets, $C$-sets 
near an idempotent. 
Dynamical characterization of central sets near an idempotent, quasi-central sets near an idempotent 
and $C$-sets near an idempotent is established in section 2, section 3 and section 4 respectively.\\

\section {dynamical characterization of central set near an idempotent}

Let us start this section with the following  definition of central 
set near an idempotent which was introduced by M. A. Tootkaboni, T. Vahed in (\cite[Definition 3.9(a)]{mat}).\\

\textbf{Definition 2.1.} Let $S$ be a dense subsemigroup of $(T,+)$ and $e\in E(T)$. Let  $\tau_e(S)=\{U\cap S:U\in \tau_e\}$.
A set $A\subseteq S$ is central set near $e$ if and only if there is some 
idempotent $p\in K$ with $A \in p$.\\

Following definition is \cite[Definition 3.3(b)]{mat}.\\

\textbf{Definition 2.2.} Let $S$ be a dense subsemigroup of $(T, +)$ and $e\in E(T)$. A subset 
$B$ of $S$ is syndetic near $e$ if and only if for each $U\in \tau_e$, there exist
some $F \in P_f(U\cap S)$ and some $V\in \tau_e$ such that 
$S\cap V\subseteq \bigcup_{t \in F}(-t+B)$.\\

S. K. Patra defined uniform recurrence and proximality near zero in \cite[Definition 2.3(a), (b)]{sou}.\\

\textbf{Definition 2.3.} Let S be a dense subsemigroup of $((0,\infty),+)$ and 
$(X,\langle T_s\rangle _{s\in S})$ be a dynamical system.

(a) A point $x \in X$ is a uniformly recurrent point near zero if and only if 
for each neighbourhood $W$ of $x$, $\{s \in S: T_s \in W  \}$ is syndetic near zero.

(b) Points $x$ and $y$ of $X$ are proximal near zero if and only if for every 
neighbourhood $U$ of the diagonal in $X\times X$ and for each $\epsilon>0$ there 
exists $s\in S \cap (0, \epsilon)$ such that $(T_s(x), T_s(y)) \in U$.\\

We now state \cite[Theorem 19.11]{hindalg}.\\

\textbf{Theorem 2.4.} Let $(X,\langle T_s \rangle_{s\in S})$ be a dynamical system and define 
$\theta : S \rightarrow X^X$ by $\theta(s)=T_s$. Then $\tilde{\theta}$ is a 
continuous homomorphism from $\beta S$ onto the enveloping semigroup of 
$(X,\langle T_s\rangle _{s\in S})$. ($\tilde{\theta}$ is the continuous extension of $\theta$).\\

The following notation will be convenient in the next section.\\

\textbf{Definition 2.5.} (\cite[Definition 19.12]{hindalg}) Let $(X,\langle T_s\rangle _{s\in S})$ be a dynamical 
system and define $\theta : S\rightarrow X^X$ by $\theta(s)=T_s$. For each 
$p \in \beta S$, let $T_p= \tilde{\theta}(p)$.\\

As an immediate consequence of Theorem 2.4 we have the following remark  \cite[Remark 19.13]{hindalg}.\\

\textbf{Remark 2.6.} Let $(X,\langle T_s\rangle _{s\in S})$ be a dynamical system and let 
$p,q \in \beta S$. Then $T_p\circ T_q=T_{pq}$ and for each 
$x \in X$, $T_p(x)$=$p$-$\lim_{s \in S}T_s(x)$.\\

Now we shall introduce the notion of uniform recurrence and proximality near $e$.\\

\textbf{Definition 2.7.} Let S be a dense subsemigroup of $(T,+)$ and $e\in E(T)$. Let 
$(X,\langle T_s\rangle _{s\in S})$ be a dynamical system.

(a) A point $x \in X$ is a uniformly recurrent point near $e$ if and only if 
for each neighbourhood $W$ of $x$, $\{s \in S: T_s(x) \in W  \}$ is syndetic near $e$.

(b) Points $x$ and $y$ of $X$ are proximal near $e$ if and only if for every 
neighbourhood $U$ of the diagonal in $X\times X$ and for each $V\in \tau_e$ there 
exists $s\in S \cap V$ such that $(T_s(x), T_s(y)) \in U$.\\

\textbf{Lemma 2.8.} Let $S$ be a dense subsemigroup of $(T,+)$ and $e\in E(T)$. 
Let $(X,\langle T_s\rangle _{s\in S})$ be a dynamical system and let 
$x, y\in X$. Then $x$ and $y$ are proximal near $e$ if and only if there 
is some $p \in e^{*}_S$ such that $T_p(x)=T_p(y)$.\\ 

\textbf{Proof.} Similar proof as the proof of \cite[Lemma 2.7]{sou}.\\

As a consequence of the above lemma one can say that, points $x$ and $y$ of $X$ 
are proximal near $e$ if and only if there is some $p \in e^{*}_S$ such that $T_p(x)=T_p(y)$.\\

\textbf{Lemma 2.9.} Let $S$ be a dense subsemigroup of $(T,+)$ and $e\in E(T)$. 
Let $(X,\langle T_s\rangle _{s\in S})$ be a dynamical system and $L$ be a minimal left ideal of $e^{*}_S$ and $x \in X$.

The following statements are equivalent :

(a) The point $x$ is a uniformly recurrent point near $e$ of $(X,\langle T_s\rangle _{s\in S})$.

(b) There exists $u \in L$ such that $T_u(x)=x$.

(c) There exist $y \in X$ and an idempotent $u \in L$ such that $T_u(y)=x$.

(d) There exists an idempotent $u \in L$ such that $T_u(x)=x$.\\

\textbf{Proof.} (a)$\Rightarrow$(b) Pick any $v \in L$. For each $U\in \tau_x$ and $x$ in $X$, let $B_U = \{s \in S : T_s(x) \in U\}$. 
Each $B_U$ is syndetic near near $e$ as $x$ is a uniformly recurrent point near $e$. Therefore, for each $W \in \tau_e$ there are some 
$F_{(U, W)} \in P_f(W \cap S)$ and some $V\in \tau_e$ 
such that $S \cap V \subseteq \bigcup_{t\in F_{(U, W)}}(-t + B_U )$. So for each
$U \in \tau_x$ and each $W \in \tau_e$, pick $t_{(U,W)} \in F_{(U, W)}$ such that $-t_{(U, W)} + B_U \in v$.
For any $U \in \tau_x$ and $W \in \tau_e$, we set $C_{(U,W)} = \{t_{(V, W)} : V \in \tau_x$, $V \subseteq U \}$ and 
$C_U=\bigcup_{W\in \tau_e}C_{(U, W)}$. 
Then $\{C_U : U \in \tau_x \} \cup \{S \cap W : W \in \tau_e \}$ has the finite intersection property 
so pick $q\in e^{*}_S$ such that
$\{C_U : U \in \tau_x \} \subseteq q$ and let $u = q + v$. Since $L$ is a left ideal of $e^{*}_S$, $u \in L$. 
Now we show that $T_u(x) = x$. Let $U \in \tau_x$. If we can show that $B_U \in u$, then we are done.
If possible, let $B_U\not \in u$. Then $\{t \in S : -t + B_U \not \in v\}$ and $C_U \in q$ and so pick
$t \in C_U$ such that $-t + B_U\not \in v$. Pick $V\in \tau_x$ with $V \subseteq U$ such that $t= t_{(V, W)}$ for some $W \in \tau_e$. 
Then $-t + B_V \in v$ and $-t + B_V \subseteq −t + B_U$, a contradiction.\\

For (b)$ \Rightarrow (c)$ and (c)$\Rightarrow$(d) see the proof \cite[Lemma 2.8]{sou}.\\

(d)$\Rightarrow$(a) Let $U\in \tau_x$ and $x\in X$ and let $B=\{ s \in S:T_s(x)\in U\}$. If possible, let 
 $B$ is not syndetic near $e$. Then there exists $W \in \tau_e$ such that 
$\{ S\setminus \bigcup_{t \in F}(-t+B)$: $F \in P_f(W \cap S)\}$ $\cup $ $\{S \cap V : V \subseteq T \}$
has the finite intersection property. So pick some $w \in e^{*}_S$ such that 
$\{ S\setminus \bigcup_{t \in F}(-t+B)$: $F \in P_f(W \cap S)\}$ $\subseteq w$.\\
We claim that $(e^{*}_S+w)\cap ClB=\emptyset$. If not, then there are some 
$v \in e^{*}_S$ with $B \in v+w$. Then one can get some $t \in P_f(S \cap W)$ 
with $-t+B \in w$, a contradiction. Therefore $(e^{*}_S+w)\cap ClB=\emptyset$.\\ 
Let $L^{'}= e^{*}_S+w$. Then $L^{'}$ is a left ideal 
of $e^{*}_S$, so $L^{'}+u$ is a 
left ideal of $e^{*}_S$ which is contained in $L$, and hence by minimality of $L$, we have $L^{'}+u=L$. Thus we 
can pick some $v \in L^{'}$ such that $v+u=u$. By using Remark 2.6, 
$T_v(x)=T_v(T_u(x))=T_{v+u}(x)=T_u(x)=x$, so in particular $B \in v$. But, $v \in L^{'}$ 
and $L^{'}\cap ClB=\emptyset$, a contradiction.\\

\textbf{Lemma 2.10.} Let $S$ be a dense subsemigroup of $(T,+)$ and $e \in E(T)$. 
Let $(X,\langle T_s\rangle _{s \in S})$ be a dynamical system and let $x \in X$. 
Then for each $W \in \tau_e$ there is a uniformly recurrent point $y$ near $e$ with
$y \in Cl\{ T_s(x): s \in S \cap W\}$ such that $x$ and $y$ are 
proximal near $e$.\\

\textbf{Proof.} Let $L$ be any minimal left ideal of $e^{*}_S$ and pick an 
idempotent $u \in L$. Let $y=T_u(x)$. For each $W \in \tau_e$, clearly  
$y \in Cl\{ T_s(x): s \in S \cap W\}$. By Lemma 2.9, $y$ is a 
uniformly recurrent point near $e$ of $(X,\langle T_s\rangle _{s \in S})$. 
By Remark 2.6 we have $T_u(y)=T_u(T_u(x))=T_{u+u}(x)=T_u(x)$ . So by Lemma 2.8 $x$ and $y$ are proximal near $e$ .\\

\textbf{Lemma 2.11.} Let $S$ be a dense subsemigroup of $(T,+)$ and $e \in E(T)$. 
Let $(X,\langle T_s\rangle _{s \in S})$ be a dynamical system and let $x,y \in X$. 
If $x$ and $y$ are proximal near $e$, then there is a minimal left ideal 
$L$ of $e^{*}_S$ such that $T_u(x)=T_u(y)$ for all $u \in L$.\\

\textbf{Proof.} Choose $v \in e^{*}_S$ such that $T_v(x)=T_v(y)$ and pick a 
minimal left ideal $L$ of $e^{*}_S$ such that $L \subseteq e^{*}_S+v$. We want to show that $L$ is as required.
Let $u \in L$ and choose $w \in e^{*}_S$ such 
that $u=w+v$. Then again using Remark 2.6, we have 
$T_u(x)=T_{w+v}(x)=T_w(T_v(x))=T_w(T_v(y))=T_{w+v}(y)=T_u(y)$.\\

\textbf{Lemma 2.12.} Let $S$ be a dense subsemigroup of $(T,+)$ and $e \in E(T)$. 
Let $(X,\langle T_s\rangle_{s \in S})$ be a dynamical system and let $x,y \in X$. There is an idempotent $u$ in 
$K(e^{*}_S)$ such that $T_u(x)=y$ if and only if both $y$ is uniformly recurrent 
near $e$ and $x$ and $y$ are proximal near $e$.\\

\textbf{Proof.} 
$(\Rightarrow)$
Since $u$ is a minimal idempotent of $e^{*}_S$, there is a minimal left ideal 
$L$ of $e^{*}_S$ such that $u \in L$. Thus by Lemma 2.9, $y$ is uniformly 
recurrent near $e$. By Remark 2.6, 
$T_u(y)=T_u(T_u(x))=T_{u+u}(x)=T_u(x)$. So $x$ and $y$ are proximal near $e$\\
Converse part: Pick by Lemma 2.11 a minimal ideal $L$ of $e^{*}_S$ such that 
$T_u(x)=T_u(y)$ for all $u \in L$. Pick by Lemma 2.9 an idempotent $u \in L$ 
such that $T_u(y)=y$.\\

 The following theorem describes a dynamical characterization of Central sets near $e$. \\ 
\textbf{Theorem 2.13.} Let $S$ be a dense subsemigroup of $(T,+)$ and $e \in E(T)$ also let 
$B \subseteq S$. Then $B$ is central near $e$ if and only if there exists a 
dynamical system $(X,\langle T_s\rangle {s \in S})$ and there exist $x,y \in X$ 
and a neighbourhood $U$ of $y$ such that $x$ and $y$ are proximal near $e$, $y$ 
is uniformly recurrent near $e$, and $B=\{ s \in S: T_s(x) \in U   \}$.\\

\textbf{Proof.} $(\Rightarrow)$ Let $Q=S\cup \{e\}$, $X=Q_{\{0,1\}}$ and 
for $s \in S$ we define $T_s:X\rightarrow X$ by $T_s(x)(t)=x(t+s)$ for all $t \in Q$.
By \cite[Lemma 19.14]{hindalg}, $T_s$ is continuous. Now let $x=\chi_{B}$, the characteristic 
function of $B \subseteq Q$. That is, $x(t)=1$ if and only if $t \in B$. Choose a minimal idempotent $u$
in $e^{*}_S$ such that $B \in u$ and let $y=T_u(x)$. Then by Lemma 2.12, $y$ is 
uniformly recurrent near $e$ and $x$ and $y$ are proximal near $e$.\\
Now let $U=\{ z \in X:z(e)=y(e) \}$. Then $U$ is a neighbourhood of $y$ in $X$. 
We note that $y(e)=1$. Indeed, $y=T_u(x)$ so, $\{ s \in S:T_s(x) \in U  \}\in u$ 
and we may choose some $s \in B$ such that $T_s(x) \in U$. Then $y(e)=T_s(x)(e)=x(s+e)= x(s)=1$.
Thus given any $s \in S$, $s \in B\Leftrightarrow x(s)=1 \Leftrightarrow T_s(x)(0)=1 \Leftrightarrow T_s(x)\in U$.
\\
$(\Leftarrow)$ Choose a dynamical system $(X,\langle T_s\rangle {s \in S})$, points $x,y \in X$
and a neighbourhood $U$ of $y$ such that $x$ and $y$ are proximal near $e$ with 
$y$ uniformly recurrent near $e$ and $B= \{ s \in S : T_s(x) \in U  \}$. 
Choose by Lemma 2.12 a minimal idempotent $u$ in $e^{*}_S$ such that $T_u(x)=y$. Then $B \in u$.\\

\section{Dynamical characterization of quasi-central sets near an idempotent}
In this present section we shall define quasi-central sets near idempotent and deduce its dynamical characterization near idempotent. \\

\textbf{Definition 3.1.} Let $S$ be a dense subsemigroup of $(T, +)$ and  $e \in E(T)$. 
Then $C$ is said to be quasi central near $e$ if and only if there is an idempotent 
$p$ in $Cl\text{ }K(e^{*}_S)$ such that $C \in p$.\\

 We need the following two definitions to characterize a quasi-central sets near an idempotent.\\

\textbf{Definition 3.2.} (\cite[Definition 3.5]{mat})
Let $S$ be a dense subsemigroup of $(T, +)$ and  $e \in E(T)$. A subset $B$ of $S$ is topologically
piecewise syndetic near $e$ if and only if for each $U \in \tau_e$ there exist some $F \in P_f(U\cap S)$ and some $V\in \tau_e$ 
such that for each $G \in P_f( S)$ and $O\in \tau_e$ there exists $x \in S \cap O$ 
such that 
$(G \cap V) + x \subseteq \bigcup_{t \in F}(-t+B)$.

Let us now introduce the notion of jointly interminittently uniform recurrence near $e$.\\

\textbf{Definition 3.3.} Let $(X, \langle T_s\rangle_{s\in S})$ be a dynamical system 
and let $x,y \in X$. The pair $(x,y)$ is jointly interminittently uniformly recurrent 
near $e$ (abbreviated as $JIUR_e$) if and only if for every neighbourhood $U$ of $y$, 
the set $\{s \in S : T_s(x) \in U \text{ and } T_s(y)\in U\}$ is topologically piecewise syndetic near $e$.\\

We state the following theorem which serves our purpose.\\

\textbf{Theorem 3.4.}(\cite[Theorem 3.6]{mat})Let $S$ be a dense subsemigroup of $(T, +)$ and  $e \in E(T)$ and  $A \subseteq S$. 
Then $K \cap Cl_{\beta S }(A) \neq \emptyset$ if and only if $A$ is topologically piecewise syndetic near $e$.\\

\textbf{Lemma 3.5.}
Let $S$ be a dense subsemigroup of $(T, +)$ and  $e \in E(T)$ and 
let $\mathcal{K} = \{A \subseteq S : S \setminus A \text{ is not topologically piecewise syndetic near $e$}\}$. 
Then $\mathcal{K}$ is a filter on $S$ with $Cl \text{ } K(e^{*}_S) = \overline{\mathcal{K}}$, 
which is a compact subsemigroup of $\beta S$.\\

\textbf{Proof.} Since $\mathcal{K} = \bigcap K(e^{*}_S)$, by \cite[Theorem 3.20(b)]{hindalg} we have   
$\mathcal{K}$ is a filter and $\overline{\mathcal{K}} = ClK(e^{*}_S)$. By \cite[Theorem 2.15]{hindalg} $ClK(e^{*}_S)$ is a right ideal of $e^{*}_S$, 
so in particular $\overline{\mathcal{K}}$ is a compact subsemigroup of $e^{*}_S$. Therefore
$Cl \text{ } K(e^{*}_S) = \overline{\mathcal{K}}$ 
 is a compact subsemigroup of $\beta S$. \\

In the following theorem we shall give a dynamical characterization of quasi-central sets near $e$.\\

\textbf{Theorem 3.6.} Let $S$ be a dense subsemigroup of $(T, +)$ and  $e \in E(T)$ and let $A \subseteq S$. 
The set $A$ is quasi-central near $e$ if and only if there exist a dynamical system $(X,\langle T_s \rangle_{s \in S})$, 
points $x$ and $y$ in $X$, and a neighbourhood $U$ of $y$ such that the the
pair $(x,y)$ is $JIUR_e$ and $A = \{s \in S : T_s(x) \in U\}$.\\

\textbf{Proof.}\\ Let $\mathcal{K} = \{B \subseteq S : S\setminus B \text{ is not a topologically piecewise syndetic set near $e$}\}$ 
and note that $\mathcal{L}(\mathcal{K}) = \{A \subseteq S : A \text{ is topologically piecewise syndetic near }$e$\}$. By Lemma 3.5, 
we have $\mathcal{K}$ is a filter and $\overline{\mathcal{K}} = Cl \text{ } K(e^{*}_S)$ which is a compact subsemigroup of $\beta S$. 
Now we can apply Theorem 1.8 to prove our required statement.\\

\section{Dynamical characterization of $C$-sets near an idempotent}
For the discussion of $C$-sets near an idempotent we need some definitions, lemma and theorems.\\

\textbf{Definition 4.1.}(\cite[Definition 3.1]{mat})
Let $(T, +)$ be a semitopological semigroup.\\
(a)Let $\mathcal{B}$ be a local base at a point $x\in T$. We say $\mathcal{B}$ has the finite cover property if 
$\{V\in \mathcal{B}: y\in V\}$ is finite for each $y\in T-\{x\}$.\\
(b)Let $S$ be a dense subsemigroup of $T$ and $e\in E(T)$. Let $\{x_n \}$ be a sequence in $S$. We say $\sum_{n=1}^{\infty}x_n$ converges near $e$ 
if for each $U\in \tau_e$ there 
exists $m\in \mathbb{N}$ such that $FS(\langle x_n \rangle_{n=k}^l) \subseteq U$ for each $l>k\geqslant m$.\\
(c)Let $\mathcal{B}=\{U_n: n\in \mathbb{N}\}$ be a countable local base  at the point $x\in T$ such that $U_{n+1}\subseteq U_n$ for each $n\in \mathbb{N}$,
$U_{n+1}+U_{n+1}\subseteq U_n$ for each $n\in \mathbb{N}$, and for each sequence $\{x_n \}$ if $x_n \in U_n$ for each $n\in \mathbb{N}$ then $\sum_{n=1}^{\infty}x_n$ converges near $x$.
Then we say $\mathcal{B}$ is a countable local base for convergence at the point $x\in T$.\\
(d)Let $\mathcal{B}$ be a local base at the point $x\in T$. If $\mathcal{B}$ satisfies in conditions (a) and (c) then $\mathcal{B}$ is called a countable local base that has the 
finite cover property for convergence at the point $x$. For simplicity we say $\mathcal{B}$ has the $\mathbf{F}$ property at the point $x$.\\

From now on we take $T$ as a commutative semigroup.\\

\textbf{Definition 4.2.}(\cite[Definition 4.1]{mat})Let $(T, +)$ be a semitopological semigroup and $S$ be a dense subsemigroup of $T$.
Let $\mathcal{B}=\{U_n \}_{n=1}^{\infty}$ has the $\mathbf{F}$ property at the point $e\in E(T)$.\\
(a)$\Phi=\{f :\mathbb{N}\rightarrow \mathbb{N}:$ for all $ n\in \mathbb{N}, f(n)\le n\}$.\\
(b)$\mathcal{Y}=\{\langle\langle{y_{i,t}}\rangle_{t=1}^{\infty}\rangle_{i=1}^{\infty}:$ for each 
$i\in \mathbb{N},\langle{y_{i,t}}\rangle_{t=1}^{\infty}$ is a sequence in $S$ and  $\sum_{t=1}^{\infty}y_{i,t}$ converges$\}$.\\
(c)Given $Y=\langle\langle{y_{i,t}}\rangle_{t=1}^{\infty}\rangle_{i=1}^{\infty}$ in $\mathcal{Y}$ and $A\subseteq S$, $A$ is a $J_Y$-set 
near $e$ if and only if for all $ n\in \mathbb{N}$ there exist $a\in U_n$ and $H\in P_f(\mathbb{N})$ such that min $H\geq n$ and for each 
$i\in \{\{1,2, \cdots ,n\}$, ${a+\sum_{t\in H}y_{i,t}}\in A$.\\
(d)Given $Y\in \mathcal{Y}$, $J_Y=\{p\in e^{*}_S:$ for all $A\in p, A$ is a $J_Y$-set near $e\}$.\\
(e)$J=\bigcap_{Y\in \mathcal{Y}}J_Y$.\\

We now state the following lemma\cite[Definition 4.2]{mat}.\\
\textbf{Lemma 4.3.}Let $(T, +)$ be a semitopological semigroup and $S$ be a dense subsemigroup of $T$.
Let $\mathcal{B}=\{U_n \}_{n=1}^{\infty}$ has the $\mathbf{F}$ property at the point $e\in E(T)$. Let $Y\in \mathcal{Y}$. Then $K\subseteq J_Y$.\\

Now we present the central sets theorem near $e$.\\

\textbf{Theorem 4.4.}Let $(T, +)$ be a semitopological semigroup and $S$ be a dense subsemigroup of $T$.
Let $\mathcal{B}=\{U_n \}_{n=1}^{\infty}$ has the $\mathbf{F}$ property at the point $e\in E(T)$ and $A$ be a central set near $e$. Suppose that 
$Y=\langle\langle{y_{i,t}}\rangle_{t=1}^{\infty}\rangle_{i=1}^{\infty}\in \mathcal{Y}$. Then there exist sequences 
$\langle{a_n}\rangle_{n=1}^{\infty}$
and $\langle {H_n}\rangle_{n=1}^{\infty}$ in $P_f(\mathbb{N}$ such that\\
(a)for each $n\in \mathbb{N}$, $a_n\in U_n$ and max $H_n<$min $H_{n+1}$, and\\
(b)for each $f\in \Phi$, $FS(\langle{a_n}+\sum_{t\in H_n}y_{f(n),t}\rangle_{n=1}^{\infty})\subseteq A$.\\

\textbf{Proof.} See Theorem 4.3 in \cite{mat}.\\

Like discrete case we now define a $C$-set near an idempotent which satisfy the conclusion of central sets theorem near an idempotent.

We start by giving the combinatorial definition of $C$-set near an idempotent. As the following combinatorial definition is rather complicated, 
we shall soon state an algebraic characterization by showing that $C$-sets near an idempotent are members of idempotents 
in a certain compact subsemigroup.\\
From now on $\mathcal{T}_e$ deontes the set of sequences in $S$ converging to $e$.

\textbf{Definition 4.5.}Let $(T, +)$ be a semitopological semigroup and $e\in E(T)$.
Let $S$ be a dense subsemigroup of $(T, +)$ and let $A \subseteq S$.
We say that $A$ is a $C$-set near $e$ if and only if for each $U\in \tau_e$, there exist functions $\alpha_U : P_f(\mathcal{T}_e) \rightarrow S$ and 
$H_U : P_f(\mathcal{T}_e) \rightarrow P_f(\mathbb{N})$ such that\\
(1) $\alpha_U(F)\in U$ for each $F \in P_f(\mathcal{T}_e)$,\\
(2) If $F, G \in P_f(\mathcal{T}_e)$ and $F \subsetneq G$, then $\max H_U(F) < \min H_U(G)$ and\\ 
(3) Whenever $m \in \mathbb{N}$, $G_1,G_2, \cdots , G_m \in P_f(\mathcal{T}_e)$, $G_1 \subsetneq G_2\subsetneq \cdots \subsetneq G_m$ and for each $i \in \{1,2, \cdots ,m\}$, $f_i \in G_i$, 
one has $$\sum_{i=1}^{m}\big(\alpha_U(G_i)+ \sum_{t \in H_U(G_i)}f_i(t)\big) \in A.$$\\

\textbf{Definition 4.6.}Let $(T, +)$ be a semitopological semigroup and $e\in E(T)$.
Let $S$ be a dense subsemigroup of $(T, +)$ and let $A \subseteq S$.\\
(1) $A$ is said to be a \textit{J}-set near $e$ if and only if whenever $F \in P_f(\mathcal{T}_e)$ and for each $U\in \tau_e$, 
there exist $a\in S \cap U$ and $H \in P_f(\mathbb{N})$ such that for each $f \in F$, $a + \sum_{t \in H}f(t) \in A$.\\
(2) $J_e(S) = \{p \in e^{*}_S : \text{ for all } A \in p, A \text{ is a }J\text{-set near $e$}\}$.\\

\textbf{Lemma 4.7.}Let $(T, +)$ be a semitopological semigroup and $e\in E(T)$.
Let $S$ be a dense subsemigroup of $(T, +)$ and let $A_1$, $A_2$ be subsets of $S$. 
If $A_1 \cup A_2$ is a $J$-set near $e$ then either $A_1$ or $A_2$ is a $J$-set near $e$.\\

\textbf{Proof.} The proof is similar to \cite[Lemma $3.8$]{bayat}.\\

\textbf{Theorem 4.8.}Let $(T, +)$ be a semitopological semigroup and $e\in E(T)$.
Let $S$ be a dense subsemigroup of $(T, +)$. Then $J_e(S)$ is a compact subsemigroup of $\beta S$.\\

\textbf{Proof.}The proof is similar to \cite[Lemma $3.9$]{bayat}.\\

\textbf{Definition 4.9.} Let $(X, \langle T_s\rangle_{s \in S})$ be a dynamical system and $x,y \in X$. 
The pair $(x,y)$ is jointly almost uniformly recurrent near $e$ (abbreviated $JIAUR_e$) if and only if for every neighbourhood $U$ of $y$, 
$\{s \in S : T_s(x) \in U \text{ and } T_s(y) \in U\}$ is a $J$-set near $e$.\\

\textbf{Lemma 4.10.}Let $(T, +)$ be a semitopological semigroup and $e\in E(T)$.
Let $S$ be a dense subsemigroup of $(T, +)$ and $\mathcal{K}= \{A \subseteq S : S \setminus A \text{ is not a }J\text{-set near $e$}\}$. 
Then $\mathcal{K}$ is a filter on $S$ with $J_e(S) = \overline{\mathcal{K}}$.\\

\textbf{Proof.} Notice that $\mathcal{K}$ is non-empty. It is also clear that $\mathcal{K}$
does not contain the empty set and is closed under super sets. Using Lemma 4.7, we say that $\mathcal{K}$ is closed under finite intersection.

Under the assumption that $\mathcal{K}$ is a filter, we have 
$\mathcal{L}(\mathcal{K}) = \{ A \subseteq S : A \text{ is a }J\text{-set} \text{near  $e$}\}$. 
By the Theorem $1.6$ we have  $J_e(S) = \overline{\mathcal{K}}$.\\

The following lemma characterizes $C$-sets near $e$ in terms of idempotents in $J_e(S)$.\\

\textbf{Lemma 4.11.} Let $(T, +)$ be a semitopological semigroup and $e\in E(T)$. Let $S$ be a dense subsemigroup of 
$(T, +)$ and $A \subseteq S$. Then $A$ is a $C$-set near $e$ if and only if there exists an idempotent $p \in J_e(S)$ such that $A \in p$.\\

\textbf{Proof.} See Theorem 3.14 in \cite{bayat}.\\

Now we give a dynamical characterization of $C$-set near $e$ in the following theorem.\\

\textbf{Theorem 4.12.}Let $(T, +)$ be a semitopological semigroup and $e\in E(T)$. Let $S$ be a dense subsemigroup of $(T, +)$ 
and let $A \subseteq S$. 
The set $A$ is $C$-set near $e$ if and only if there exist a dynamical system $(X, \langle T_s \rangle_{s \in S})$, 
points $x$ and $y$ in $X$, and a neighbourhood $U$ of $y$ such that the pair $(x,y)$is $JIAUR_e$  and $A = \{s \in S : T_s(x) \in U\}$.\\

\textbf{Proof.} Let $\mathcal{K} = \{B \subseteq S : S \setminus B \text{ is not a }J \text{-set near $e$}\}$. 
 Lemma 4.11 characterizes $C$-sets near $e$ in terms of idempotents in $J_e(S)$. Again by Lemma 4.10 we have $J_e(S)= \overline{\mathcal{K}}$.
 Now we can apply Theorem 1.8 to prove our statement.\\

\textbf {Acknowledgment.}  The authors are thankful to Prof. Dibyendu De of University of Kalyani for his guidance,
a number of valuable suggestions, sharing thoughts for betterment of the paper.

\end{document}